\documentclass[oneside]{article}
\usepackage[utf8]{inputenc}
\usepackage[english]{babel}

\usepackage{amsthm}
\newtheorem{theorem}{Theorem}[section]
\newtheorem{lemma}{Lemma}[section]

\newtheorem{assumption}{Assumption}[section]

\newtheorem{exam}{Example}[section]

\theoremstyle{definition}

\usepackage{amsmath}
\usepackage{amsfonts}
\usepackage{amssymb}
\usepackage{graphicx}
\numberwithin{equation}{section}
\usepackage{lipsum}
\usepackage{mathtools}
\usepackage{scalefnt}
\usepackage{empheq}
\usepackage[a4paper,left=3.5cm,right=3cm,top=3cm,bottom=3cm]{geometry}
\usepackage{makeidx}
\usepackage{bm}
\makeindex

\usepackage{indentfirst}

\usepackage{fullpage}
\usepackage{amsmath}
\usepackage{amssymb}
\usepackage{graphicx}
\usepackage{subfigure}
\usepackage{color}
\usepackage{epsfig}
\usepackage{epstopdf}
\usepackage[breaklinks]{hyperref}
\usepackage{algorithm,algcompatible}
\usepackage{authblk}
\usepackage{comment}

\date{}

	\title{Adapted Lie splitting method for convection-diffusion problems with singular convective term}
\author[1]{Thi Tam Dang\thanks{Thi-Tam.Dang@uibk.ac.at}}
\author[1]{Trung Hau Hoang\thanks{Trung-Hau.Hoang@uibk.ac.at}}
\author[2]{Giandomenico Orlandi\thanks{giandomenico.orlandi@univr.it}}
\affil[1]{Department of Mathematics, University of Innsbruck, Austria}
\affil[2]{Dipartimento di Informatica, Università di Verona, Italy}

\begin{document}
	\maketitle

\begin{abstract}
Splitting methods are a widely used numerical scheme for solving convection-diffusion problems. However, they may lose stability in some situations, particularly when applied to convection-diffusion problems in the presence of an unbounded convective term. In this paper, we propose a new splitting method, called the ``Adapted Lie splitting method'', which successfully overcomes the observed instability in certain cases. Assuming that the unbounded coefficient belongs to a suitable Lorentz space, we show that the adapted Lie splitting converges to first-order under the analytic semigroup framework. Furthermore, we provide numerical experiments to illustrate our newly proposed splitting approach.
\end{abstract}

\section{Introduction}
Let $\Omega $ be an open bounded subset of $\mathbb{R}^{N}$ with $N\ge 2$. We consider convection-diffusion problems of the form
\begin{equation}\label{1.1} \left\{
	\begin{aligned}
	&	\partial_{t}u(t,x) -D u(t,x) - \nabla \cdot ( c(x) u(t,x))  = f(t,x), \ \ (t,x) \in [0,T]\times \Omega,\\
	& u(t,x)|_{\partial \Omega} = b(t,x),\\
	& u(0,x) = u_0(x),	
	\end{aligned}
	\right.
\end{equation}
where $D$ denotes a second-order strongly elliptic differential operator, the convective coefficient $c(x)$ is possibly unbounded and belongs to the weak-Marcinkiewicz space $ L^{N, \infty}(\Omega)$. We impose a possibly time-dependent Dirichlet boundary condition $b: [0,T] \times \partial \Omega \rightarrow \mathbb{R}$, which is assumed to be continuously differentiable. In the context of convection-diffusion problems, the periodic boundary conditions are of significance within theoretical analyses, whereas the imposition of Dirichlet boundary conditions provides a practical approach to enforcing desired behaviors or constraints on the system, consistent with the realities of a wide range of physical and engineering problems.

In the context of linear homogeneous systems, equation \eqref{1.1} is commonly referred to as the Fokker--Planck equation, which characterizes the evolution of certain Brownian motion processes (see \cite{Porretta2015}). When the convective coefficient $c(x)$ remains bounded, the investigation of problem \eqref{1.1} has garnered attention from numerous researchers. However, in various applications, particularly within the realm of semiconductor device diffusion models, the boundedness of the convective field may not be assured. The theoretical analysis of \eqref{1.1} and its generalisations has led to the introduction of several notions of solution, notably including the concept of renormalized solution (see, e.g., \cite{Cardaliaguet2012,Porretta2015}) and entropy solution (see \cite{Boccardo2003}), that are designed to address scenarios where the convective coefficient lacks boundedness, and thus provide a particularly well-suited framework for analyzing the problem under consideration also from a numerical point of view.


In the field of computational mathematics and scientific computing it is crucial to develop effective methods for solving convection-diffusion problems like \eqref{1.1}. Among the various methods being studied, splitting methods have received significant attention as a particularly promising approach (see, e.g., \cite{hairer2013geometric,HUNDSDORFER1995191,mclachlan_quispel_2002}). These methods decompose the problem into more manageable subproblems and offer a way to improve computational efficiency and effectiveness. Furthermore, splitting methods have an advantage in preserving some relevant features of the solution such as positivity, as long as the solutions of the subproblems also enjoy that feature (see \cite{HANSEN20121428}). This makes splitting methods appealing in practical applications where for instance maintaining positivity is crucial, such as modeling physical quantities like concentrations or densities. A standard splitting approach involves splitting the convection-diffusion problem \eqref{1.1} into the following two subproblems:
\begin{equation}\label{1.2}\left\{
	\begin{aligned}
		&\partial_{t} v(t,x) = D v(t,x) + f(t,x), \\
		&	v(t,x)|_{\partial \Omega}= b(t,x),\\	
	\end{aligned}
	\right.
\end{equation}	
and
\begin{equation}\label{1.3}\left\{
	\begin{aligned}
		&	\partial_{t} w(t,x)=  \nabla \cdot (c(x) w(t,x)),\\
		&	w(t,x)|_{\Gamma} = b_1(t,x),
	\end{aligned}
	\right.
\end{equation}	
where $\Gamma$ denotes the  inflow boundary and $b(t,x)=b_1(t,x)$ for $x \in \Gamma$.	
The diffusion problem \eqref{1.2} can be efficiently treated by a fast Poisson solver (see \cite{MCKENNEY1995348}) (with complexity $\mathcal{O}(n \log n)$) where $n$ is the number of degrees of freedom. The pure convection problem \eqref{1.3} can be solved by appropriate methods. Our study focuses on the Lie splitting method, which is a well-known first-order method (see, e.g., \cite{HUNDSDORFER1995191,mclachlan_quispel_2002}), for solving the convection-diffusion problem \eqref{1.1}.  Given an initial data $u_{n}$, fixing a time step size $\tau>0$, one step of the classical Lie splitting iteration in the time interval $\left[ t_{n}, t_{n+1}\right]$ is given by
\begin{align}\label{1.4}
	u_{n+1} = \mathcal{L}_{\tau}u_{n}=  \varphi_{\tau}^{D,f} \circ \varphi_{\tau}^{C}(u_{n}),
\end{align}
where $\varphi_{t}^{D,f}$ and $\varphi_{t}^{C}$ are the solution operators of problems \eqref{1.2} and \eqref{1.3} respectively. The operator $\mathcal{L}_{\tau}$ defined in \eqref{1.4} is referred to as Lie splitting with time step size $\tau$.

Splitting methods are shown to be effective when applied to convection-diffusion problems with a bounded convective term, i.e. when the coefficients are assumed to be sufficiently smooth (see \cite{HUNDSDORFER1995191}). However, in scenarios featuring unbounded convective coefficients many additional challenges are introduced, both for the theoretical analysis and from the numerical perspective. Our numerical experiments (see Section~5 below) have shown a notable instability of splitting methods when confronted with the presence of unbounded convective terms, due to singularity effects. This observation has prompted us to investigate potential modifications of splitting methods in order to minimize the instability encountered. Efforts to develop reliable numerical techniques that have the ability to handle the complexities of unbounded convective terms are motivated by the need to improve the accuracy and dependability of computational simulations in this field.

In this paper we present the ``Adapted Lie splitting method'', an innovative approach developed to tackle the instability issue commonly found in the standard Lie splitting approach \eqref{1.4}. The proposed method aims to reduces the impact of singularity effects by decomposing the unbounded convective coefficient $c(x)$ into a combination of a bounded function and a $L^{N,\infty}(\Omega)$ function for which there is a control on its norm. Our newly approach thus centers around reformulating problem \eqref{1.1} by strategically dividing the convection term into a controlled component and another that exhibits bounded convective behavior. This reformulation, which is detailed in Section~3, effectively avoids the inherent instability conditions encountered. Furthermore, our novel splitting approach offers a promising technique for achieving numerical stability and greatly enhancing accuracy in broader scenarios characterized by unbounded convective terms. 

The convergence results for the newly proposed "Adapted Lie splitting" method are obtained in our main result, Theorem \ref{Th: 4.1}, under suitable assumptions on the convective term (see Assumption \ref{as2}). These assumptions ensure stability properties and uniform boundedness of the solution. Additionally, we conduct a thorough convergence analysis to demonstrate that the adapted Lie splitting method is a first-order convergent method.

The paper is structured as follows: Section~2 provides a review of the definition and properties of Lorentz spaces. In Section~3, we introduce an adapted Lie splitting method that effectively overcomes instability and significantly improves the accuracy of the numerical method.  In Section~4  our main theoretical convergence analysis result is presented in Theorem \ref{Th: 4.1}, which establishes the first-order convergence of the adapted Lie splitting under suitable assumptions on the convective coefficients, in an analytic semigroup framework. A rigorous error analysis is carried out to explain the behavior observed in the numerical simulations. To demonstrate the effectiveness of the proposed splitting method, numerical experiments have been conducted in both two- and three-dimensional spaces, as outlined in Sections 5 and 6.

\section{Preliminary definitions and results}
Let $\Omega$ be a bounded domain of $\mathbb{R}^{N}$. For any $1<p,q<\infty$ and $h>0$, we define the distribution function of f given by
\begin{align*}
	\lambda_{g}(h) = \vert \left\lbrace x \in \Omega: \vert g(x) \vert >h \right\rbrace \vert.
\end{align*}
The Lorentz space $L^{p,q}(\Omega)$ consists of all measurable functions $g$ which satisfies the following quantity
\begin{align*}
	\| g \|_{p,q}^{q} := p \int_0^{\infty} \left[ \lambda_{g}(h)\right]^{q/p} h^{q-1}dh \le \infty.
\end{align*}
The Lorentz space $L^{p,q}(\Omega)$ becomes a Banach space when endowed with a norm $\| \cdot \|_{p,q}$. In the case where $p=q$, the Lorentz space $L^{p,p}(\Omega)$ is idenitified as the Lebesgue space $L^{p}(\Omega)$. For $q = \infty$, the class $L^{p, \infty}(\Omega)$ coincides with the Marcinkiewicz class weak-$L^{p}(\Omega)$, and thus
\begin{align*}
	\| g\|_{p, \infty}^{p}:= \sup_{h>0} h^{p}\lambda_{g}(h) < \infty.
\end{align*}
A typical example of an element of $L^{N, \infty}(\Omega)$ is the function $g(x) = M / \vert x \vert$, where $M>0$ is a constant. Further background on the Lorentz spaces can be found in  \cite{8403ffbb-9bb1-304a-8240-0ece5b0a98cb}.\\
For $1 \le q <p<r \le \infty$, in the Lorentz spaces we have the following continuous inclusions
\begin{align*}
	L^{r}(\Omega) \subset L^{p,q}(\Omega) \subset L^{p,r}(\Omega) \subset L^{p, \infty}(\Omega) \subset L^{q}(\Omega).
\end{align*}
Let $p, q \in [1, \infty]$ with $\frac{1}{p}+ \frac{1}{p^{\prime}}= 1, \frac{1}{q} + \frac{1}{q^{\prime}}=1$. For $f \in L^{p,q}(\Omega)$ and $g \in L^{p^{\prime},q^{\prime}}(\Omega)$ we have the H\"older-type inequality 
\begin{align}\label{2.1}
	\| fg\|_1 \le \| f\|_{p,q}\|g \|_{p^{\prime}, q^{\prime}}.
\end{align}

Functions in $L^{\infty}(\Omega)$ are not dense in $L^{p,\infty}(\Omega)$. For any $g \in L^{N, \infty}(\Omega)$ we define
\begin{align*}
	\mathrm{dist}_{L^{N, \infty}(\Omega)}(g, L^{\infty}(\Omega))= \inf_{\phi \in L^{\infty}(\Omega)} \| g-\phi\|_{L^{N, \infty}(\Omega)}.
\end{align*}
We have the  characterization 
\begin{align}
	\mathrm{dist}_{L^{N, \infty}(\Omega)}(g, L^{\infty}(\Omega)) = \lim_{ K \rightarrow \infty} \| g -T_{K}g\|_{N, \infty},
\end{align}
where $T_{K}$ denotes the truncation operator at level $\pm K$, i.e. $T_{K}(s)= \frac{s}{\vert s \vert}\min\left\lbrace \vert s \vert, K \right\rbrace $, for all $s \in \mathbb{R}$. For example, the distance between $g(x) = M /\vert x \vert$ in $L^{N, \infty}(\Omega)$ and $L^{\infty}(\Omega)$ is given by the following formula:
\begin{align*}
	\mathrm{dist}_{L^{N, \infty}(\Omega)}(g, L^{\infty}(\Omega)) = M \omega_{N}^{1/N}
\end{align*}
where $\omega_{N}$ denotes the Lebesgue measure of the unit ball of $\mathbb{R}^{N}$. For further details regarding the distance to $L^{\infty}(\Omega)$, we refer the reader to \cite{Carozza1997}.\\
For any given $\delta \ge 0$, we define the subset $X_{\delta}(\Omega)$ of $L^{N, \infty}(\Omega)$ consists of all functions $g \in L^{N, \infty}$ such that
\begin{align}\label{2.2}
	\| g -T_{K}g\|_{L^{N, \infty}(\Omega)} \le \delta.
\end{align}
The existence  and uniqueness of a solution to the problem \eqref{1.1} is given in the following theorem
\begin{theorem}(Existence and uniqueness,\cite{ Farroni2021,FARRONI2018397})\label{Th: 2.2}
	There exists a proper convex set $X_{\delta}(\Omega)$ of the space $L^{N, \infty}(\Omega)$ such that, if 
	\begin{align}\label{2.5}
		c \in X_{\delta}(\Omega)
	\end{align}
	then the problem \eqref{1.1} admits a unique solution $u \in L^2(0,T; W^{1,2}(\Omega)) \cap C([0,T]; L^2(\Omega))$.
\end{theorem}
This section concludes with the Sobolev embedding theorem in Lorentz spaces, which is employed in the proof of Lemma \ref{Le: 3.1} below.
\begin{theorem}(Sobolev embedding \cite{8403ffbb-9bb1-304a-8240-0ece5b0a98cb}) \label{Th: 2.1}
	Assume that $1<p<N, q \le q \le p$, then every function $u \in W_0^{1,1}(\Omega)$ verifying $\vert \nabla u \vert \in L^{p,q}(\Omega)$ actually belongs to $L^{p^{*},q}(\Omega)$, where $p^{*}:= \frac{Np}{N-p}$ is the Sobolev exponent of $p$ and 
	\begin{align}
		\| u \|_{p^{*},q} \le S_{N,p} \| \nabla u \|_{p,q}
	\end{align}
where $S_{N,p}$ is the Sobolev constant given by $S_{N,p}= \omega_{N}^{-1/N}\frac{p}{N-p}$.
\end{theorem}

\section{Description of the method}
In this section we propose a new adaptive splitting approach for problem \eqref{1.1} that avoids the instability caused by the unbounded covective term. Let us provide first an auxiliary model problem.
\subsection{Model problem}
For any fixed $K >0$ and  $c(x) \in L^{N, \infty}(\Omega)$, define
\begin{align}\label{3.1}
	\theta(x) = \frac{T_{K}(c(x))}{c(x)},
\end{align}
where $T_{K}$ is the usual truncation operator at level $\pm K$. Consider the following auxiliary parabolic problem:
\begin{equation}\label{3.2}\left\{
	\begin{aligned}
		&\partial_{t}u(t) - Du(t) -  \nabla \cdot \Big( (1- \theta(\cdot)) u(t) \Big)  = f(t)+ \nabla \cdot ( \theta(\cdot)  u(t)),\\ 
		& u(t)|_{\partial \Omega} = b(t),\\
		& u(0)= u_0,
	\end{aligned}
	\right.
\end{equation}
where $D$ denotes a second-order strongly elliptic differential operator and $\theta(\cdot)$ denotes the space-dependent adaptive coefficient. Henceforth, we adopt the notation $u(t)$ instead of $u(t,x)$ and $\theta(\cdot)$ instead of $\theta(x)$ for the sake of notational brevity.\\
The second-order linear elliptic differential operator $D$ is defined as follows:
\begin{equation*}
	\begin{aligned}
		Du = \sum_{i,j=1}^{N}- \partial_{i}(a_{ij}(x) \partial_{j} u)  + \sum_{i=1}^{N} \alpha_{i}(x)\partial_{i}u+ \beta(x) u,
	\end{aligned}
\end{equation*}
where the matrix-value function $(a_{i,j}(x)) \in \mathbb{R}^{N \times N}$ is assumed to be symmetric, i.e. $a_{i,j}(x)= a_{j,i}(x)$, for $1 \le i,j \le N$. Further, we assume that the coefficient $a_{ij}, \alpha_{i}, \beta$ are sufficiently smooth. \\
For any $ x \in \Omega$, $\sigma \in \mathbb{R}^{N}$, the operator $D$ of the form \eqref{3.2} is strongly elliptic if there exists a constant $\eta$ such that
\begin{align}\label{3.3}
	\sum_{i,j=1}^{N} a_{ij}(x) \sigma_{i}\sigma_{j} \ge \eta \vert \sigma \vert^2.
\end{align}

The operator \eqref{3.3} on $L^2(\Omega)$ corresponds to elliptic operator in divergence form $D$ written formally as
\begin{align*}
	Du = - \mathrm{div} (a_{ij} \nabla u ) + \alpha \cdot \nabla u +  \beta u.
\end{align*}
 This operator is associated with the bilinear form 
\begin{align*}
	a : H^1(\Omega) \times H^1(\Omega) \rightarrow \mathbb{Z}
\end{align*}
given by
\begin{align*}
	a(u,v) = \int_{\Omega} \Big( \sum_{i,j=1}^{N} a_{ij}\partial_{i}u \partial_{j} \bar{v} + \alpha_{j} \partial_{j}u \bar{v} \Big) + \beta u \bar{v} dx.
\end{align*}
Let us moreover define the operator 
\begin{align*}
	\mathcal{A}u := - \mathrm{div}(a_{ij} \nabla u )+ \alpha \cdot \nabla u - \mathrm{ div} ((1- \theta(\cdot)) \nabla u) + \beta u.
\end{align*}
The operator $\mathcal{A}$ will be associated with the bilinear form 
\begin{align*}
	\hat{a} : H^1(\Omega) \times H^1(\Omega) \rightarrow \mathbb{Z}
\end{align*}
given by
\begin{align}\label{3.3a}
		\hat{a}(u,v) = \int_{\Omega} \Big( \sum_{i,j=1}^{N} a_{ij}\partial_{i}u \partial_{j} \bar{v} + \sum_{j=1}^{N} \Big( \alpha_{j} \partial_{j}u \bar{v} + (1- \theta(\cdot))u \partial_{j} \bar{v}\Big) + \beta u \bar{v} \Big)dx.
\end{align}
\begin{lemma}\label{Le: 3.1}
	The form $\hat{a}$ defined by \eqref{3.3a} is $L^2(\Omega)$-elliptic.
\end{lemma}
\begin{proof}
Let us define $\hat{a}_1: H^1(\Omega)\times H^1(\Omega) \rightarrow \mathbb{Z}$ is given by 
\begin{align*}
	\hat{a}_1(u,v) = \tilde{a}_1(u,v) + 	\int_{\Omega} \Big( \sum_{j=1}^{N} (1- \theta(\cdot))u \partial_{j} \bar{v} \Big) dx, 
\end{align*}	
where $\tilde{a}_1(u,v): H^1(\Omega)\times H^1(\Omega) \rightarrow \mathbb{Z}$ is defined as follows
\begin{align*}
	\tilde{a}_1(u,v) = 	\int_{\Omega}\Big(  \sum_{i,j=1}^{N} a_{ij}\partial_{i}u \partial_{j} \bar{v} \Big)dx.
\end{align*}
The ellipticity condition \eqref{3.3}, the assumption \eqref{2.5} with $\delta = \frac{\eta}{2 S_{N,2}}$, and the H\"older inequality \eqref{2.1}, along with the Sobolev embedding theorem \ref{Th: 2.1}, are employed to obtain
\begin{align*}
  \vert \hat{a}_1(u)\vert &= \mathrm{Re} \tilde{a}_1(u)  + \| c- T_{K}c\|_{L^{N, \infty}(\Omega)}\| u \|_{L^2(\Omega)} \| \nabla u \|_{L^2(\Omega)} \\
 & \ge \eta \| \nabla u\|_{L^2(\Omega)}^2- S_{N,2} \| c-T_{K}c\|_{L^{N, \infty}(\Omega)} \| \nabla u\|_{L^2(\Omega)}^2\\
  & \ge \frac{\eta}{2}\| \nabla u \|_{L^2(\Omega)}^2.
\end{align*}
Thus, we have
\begin{align}
\vert \hat{a}_1(u)\vert   + \frac{\eta}{2} \| u\|_{L^2(\Omega)}^2 \ge \frac{\eta}{2} \| u\|_{H^1(\Omega)}^2,
\end{align}
which implies that $\hat{a}_1$ is $L^2(\Omega)$-elliptic.\\ 
Let us define next $\hat{a}_2: H^1(\Omega)\times H^1(\Omega) \rightarrow \mathbb{Z}$ by 
\begin{align*}
	\hat{a}_2(u,v) = \int_{\Omega}\Big(  \sum_{j=1}^{N} \alpha_{j} \partial_{j}u \bar{v} + \beta u \bar{v} \Big)dx.
\end{align*}
It is obvious that  $\hat{a}= \hat{a}_1 + \hat{a}_2$.\\
The boundedness of the coefficients $\alpha_{j}, \beta$ and the Young's inequality imply that there exists a constant $C$ such that	
\begin{equation}
	\begin{aligned}
	\vert \hat{a}_2(u)\vert& \le C \int_{\Omega} \vert \nabla u \vert \vert u \vert dx + \| \beta \|_{L^{\infty}(\Omega)} \| u\|_{L^2(\Omega)}^2\\
	& \le C \| u\|_{H^1(\Omega)} \| u\|_{L^2(\Omega)}\\
	& \le  \frac{1}{2} \Big( \eta \| u\|_{H^1(\Omega)}^2 + \frac{1}{\eta}C^2 \| u\|_{L^2(\Omega)}^2\Big).
	\end{aligned}
\end{equation}
This leads to
\begin{align*}
	\mathrm{Re}(\hat{a}_1(u) + \hat{a}_2(u)) + \Big(  \frac{\eta}{2} + \frac{C^2}{2 \eta}\Big) \| u\|_{H^1(\Omega)}^2 \ge \frac{\eta}{2} \| u \|_{H^1(\Omega)}^2.
\end{align*}
This concludes the proof. 

\end{proof}

\subsection{Adapted Lie splitting method}
In order to describe the construction of our proposed splitting method, let us first define 
\begin{align}\label{3.4}
	\hat{D}u(t) = - D u(t) - \nabla \cdot  \Big( (1- \theta(\cdot)) u(t) \Big) ,
\end{align}
where $\theta(\cdot)$ is given by \eqref{3.1}. As a consequence of Lemma \ref{Le: 3.1}, the operator $\hat{D}$ remains a second-order strongly elliptic differential operator. The problem \eqref{3.2} is then reformulated in terms of $\hat{D}$ and $\theta$ as follows:
\begin{equation}\label{3.5}\left\{
	\begin{aligned}
	&	\partial_{t} u(t) = \hat{D} u(t) + \nabla \cdot ( \theta(\cdot)  u(t))  +f(t),\\
	& u(t)|_{\partial \Omega} =b(t),\\
	& u(0) = u_0.	
 	\end{aligned}
	\right.
\end{equation}
We now proceed to the splitting of \eqref{3.4} into the following two subproblems:
\begin{equation}\label{3.6}\left\{
	\begin{aligned}
	&	\partial_{t} v(t) = \hat{D} v(t) + f(t),\\
	&	v(t)|_{\partial \Omega} = b(t),
	\end{aligned}
	\right.
\end{equation}
and 
\begin{equation}\label{3.7}\left\{
	\begin{aligned}
	&	\partial_{t}w(t) = \nabla \cdot \Big( \theta(\cdot) w(t) \Big),\\
	&  w(t)|_{\Gamma} = b_1(t).
	\end{aligned}
	\right.
\end{equation}
We remark that in the splitting procedure, we impose a Dirichlet boundary condition for the diffusion equation \eqref{3.6} while an inflow boundary condition is imposed for the pure convection equation \eqref{3.7}.

We define the iterative step of the adapted Lie splitting scheme with time step size $\tau>0$ in the time interval $\left[ t_{n}, t_{n+1}\right]$  as follows:
\begin{align}\label{3.8}
	u_{n+1} = \mathcal{L}_{\tau}^{Ad}u_{n}= \varphi_{\tau}^{\hat{D},f} \circ \varphi_{\tau}^{\hat{C}}(u_{n}),
\end{align}
where $\varphi_{t}^{\hat{D},h}$ and $\varphi_{t}^{\hat{C}}$ are the solution operators of the problems \eqref{3.6} and \eqref{3.7}, respectively. The operator $\mathcal{L}_{\tau}^{Ad}$ in \eqref{3.8} is referred to as adapted Lie splitting scheme with time step size $\tau$. \\
Each step of the adapted Lie splitting scheme is explicitly carried out in the following way:
\begin{algorithm}
	\caption{Adapted Lie splitting scheme for \eqref{1.1}}
	\begin{itemize}
		\item[$1.$] Compute the solution of \eqref{3.7} with the intial value $w_{n}(0)= u_{n}$ to obtain $w_{n}(\tau)$.
		\item[$2.$] Compute the solution of \eqref{3.6} with the initial value $v_{n}(0) = w_{n}(\tau)$ to obtain $u_{n+1} =\mathcal{L}_{\tau}^{Ad}u_{n}= v_{n}(\tau)$.
	\end{itemize}	
\end{algorithm}

\section{Convergence analysis}
In this section we derive a thorough error analysis for adapted Lie splitting applied to \eqref{1.1}. More precisely, we prove that the adapted Lie splitting method is first-order convergent in the analytic semigroup setting and under appropriate assumptions.
\subsection{Analytical framework}
The convergence analysis will be conducted within the context of analytic semigroups. To this end, we will rewrite the problem \eqref{3.5} in a way that allows for the imposition of homogeneous Dirichlet boundary conditions. This transformation is only employed for the purposes of error analysis and does not affect our numerical scheme in general.
To achieve this, we introduce a smooth function, $z(t)$, which is a solution to the following elliptic problem:
\begin{equation}\left\{
	\begin{aligned}
	&	\hat{D}z(t) = 0,\\
		&z(t)|_{\partial \Omega}=b(t),
	\end{aligned}
	\right.
\end{equation}
where $\hat{D}$ denotes a second-order strongly elliptic differential operator is given by \eqref{3.4}.  \\
Let us now define $\tilde{u}(t)= u(t)-z(t)$ which satisfies the following abstract evolution equation
\begin{equation}\label{4.1}\left\{
	\begin{aligned}
	&	\partial_{t} \tilde{u}(t) +A \tilde{u}(t) = \nabla \cdot (\theta(\cdot) \tilde{u}(t)) + f(t),\\
	& u(0)=u_0,
	\end{aligned}
\right.
\end{equation}
where the operator $A$ is defined as: $A\psi= -\hat{D}\psi$, for all $\psi \in D(A)$. For example, in the case of the operator $\hat{D}$ on $L^2(\Omega)$, it can be shown that $D(A)=H^2(\Omega) \cap H_0^1(\Omega) \subset L^2(\Omega)$. It is observed that the boundary condition for the equation \ref{4.1} is included in the domain of the operator $A$.\\
We make use of the following assumption for the operater $A$.
\begin{assumption}\label{as1}
	Let $X$ be a Banach space with norm $\Vert \cdot \Vert$. We assume that $A$ is a linear operator on $X$ and $-A$ is the infinitesimal generator of an analytic semigroup $e^{-tA}$ on $X$.
\end{assumption}
In the following we recall some properties of analytic semigroups which will be used extensively in the upcoming convergence analysis. For more background on the analytic framework, we refer to the books of Henry \cite{henry1981geometric} or Pazy \cite{pazy2012semigroups}. In particular, for any $t \in [0,T]$ there exists a constant $C$ such that
\begin{align}\label{4.3}
	\Vert e^{-tA}\Vert \le C.
\end{align}
Moreover, for any $\alpha >0$, the operator $A$ enjoys the parabolic smoothing property
\begin{align}\label{4.4}
	\| A^{\alpha} e^{-tA} \| \le \frac{C}{t^{\alpha}},
\end{align}
uniformly for $t \in (0, T]$.\\
We now proceed by splitting the problem \eqref{4.1} into two partial flows:
\begin{equation}\label{4.5}
	\begin{aligned}
		\partial_{t} \tilde{v}(t)+ A \tilde{v}(t)= f(t),
	\end{aligned}
\end{equation}
and 
\begin{equation}\label{4.6}\left\{
	\begin{aligned}
		&\partial_{t}\tilde{w}(t)= \nabla \cdot( \theta(\cdot) \tilde{w}(t)),\\
		& \tilde{w}(t)|_{\Gamma} =0.
	\end{aligned}
\right.
\end{equation}
We note that the problem \eqref{4.5} with the homeogeneous Dirichlet boundary conditions included in the domain of operator $A$.\\
The exact solution of \eqref{4.1} at the time $t_{n+1} = t_{n}+\tau$ is given by the variation-of-constants formula
\begin{equation}\label{4.7}
	\begin{aligned}
		\tilde{u}(t_{n+1})= e^{-\tau A}\tilde{u}(t_{n})+ \int_0^{\tau} e^{-(\tau-s)A} \Big( \nabla \cdot (\theta(\cdot)\tilde{u}(t_{n}+s) )  + f(t_{n}+s) \Big) ds.
	\end{aligned}
\end{equation}
This gives the exact solution of \eqref{3.5} below
\begin{equation}\label{4.8}
	\begin{aligned}
		u(t_{n+1})&= z(t_{n+1}) +  e^{-\tau A} \tilde{u}(t_{n}) + \int_0^{\tau} e^{-(\tau-s)A} f(t_{n}+s)ds \\
		& \qquad + 	\int_{0}^{\tau} e^{-(\tau-s)A}   \nabla \cdot  (\theta(\cdot) \tilde{u}(t_{n}+s))  ds.
	\end{aligned}
\end{equation}
The upcoming subsections will investigate the local error and the global error of the adapted Lie splitting method applied to \eqref{3.5}.
\subsection{Local error}
We first derive a bound for the local error for the adapted Lie splitting applied to \eqref{3.5}. For this purpose, let us consider one step of the numerical solution, starting at time $t_{n}$ with the initial value $\tilde{w}(t_{n})= \tilde{u}(t_{n})$ on the exact solution. The solution of \eqref{4.6} with full step size $\tau$ can be expressed in a Taylor series with integral remainder as follows:
\begin{align}\label{4.9}
\tilde{w}_{n}(\tau)= \tilde{w}(t_{n}+\tau) = \tilde{u}(t_{n}) + \tau \nabla \cdot (\theta(\cdot) \tilde{u}(t_{n})) + \int_0^{\tau}(\tau-s) \tilde{w}^{\prime \prime}(t_{n}+ s) ds.
\end{align}
Next we integrate \eqref{4.5} with initial data $\tilde{v}(t_{n})= \tilde{w}_{n}(\tau)$ to obtain
\begin{equation}\label{4.10}
	\begin{aligned}
	\mathcal{L}_{\tau}^{Ad}\tilde{u}(t_{n})=	\tilde{v}(t_{n}+\tau) &= e^{-\tau A} \tilde{u}(t_{n})+ \int_0^{\tau} e^{-(\tau-s)A}f(t_{n}+s) ds\\
	& \qquad + \tau e^{-\tau A} \nabla \cdot (\theta(\cdot) \tilde{u}(t_{n}))+ \int_0^{\tau} e^{-\tau A} (\tau-s) \tilde{w}^{\prime \prime}(t_{n}+s)ds.
	\end{aligned}
\end{equation}
Therefore, the Lie splitting scheme applied to \eqref{3.5} in the following gives the numerical solution
\begin{equation}\label{4.11}
	\begin{aligned}
	u_{n+1}= \mathcal{L}_{\tau}^{Ad} u(t_{n}) &=  z(t_{n+1})+  e^{-\tau A} \tilde{u}(t_{n})+ \int_0^{\tau} e^{-(\tau-s)A}f(t_{n}+s) ds\\
	& \qquad +  \tau e^{-\tau A} \nabla \cdot (\theta(\cdot) \tilde{u}(t_{n}))+ \int_0^{\tau} e^{-\tau A} (\tau-s) \tilde{w}^{\prime \prime}(t_{n}+s)ds.
	\end{aligned}
\end{equation}
The local error $d_{n+1} = \mathcal{L}_{\tau}^{Ad}(t_{n})- u(t_{n+1})$ is obtained by subtracting the exact solution \eqref{4.8} from the numerical solution \eqref{4.11}. This gives the following representation of the local error:
\begin{equation}
	\begin{aligned}
		d_{n+1}& = \tau e^{-\tau A} \nabla \cdot (\theta(\cdot) \tilde{u}(t_{n})) - \int_0^{\tau} e^{-(\tau-s)A} \nabla \cdot (\theta(\cdot) \tilde{u}(t_{n}+s)) ds\\
		& \qquad + \int_0^{\tau}e^{-\tau A} (\tau -s) \tilde{w}^{\prime \prime}(t_{n}+s) ds.
	\end{aligned}
\end{equation}
We set
\begin{align}
 l_{n}(s) = e^{-(\tau-s)A}  \tilde{l}_{n}(s) , \quad  \tilde{l}_{n}(s)= \nabla \cdot (\theta(\cdot) \tilde{u}(t_{n}+s)).
\end{align}
The Taylor expansion of $l_{n}(s)$ gives
\begin{align*}
	\int_0^{\tau} l_{n}(s)ds = \tau l_{n}(0) + \int_0^{\tau} \int_0^{\sigma} l_{n}^{\prime}(\sigma) d\sigma ds.
\end{align*}
This implies that
\begin{equation}\label{4.13}
	\begin{aligned}
		d_{n+1} = \int_0^{\tau} e^{-\tau A} (\tau -s) \tilde{w}^{\prime \prime}(t_{n}+s) ds - \int_0^{\tau} \int_0^{\sigma} l_{n}^{\prime}(\sigma) d\sigma ds.
	\end{aligned}
\end{equation}
As a consequence of Theorem \ref{Th: 2.2}, the exact solution $u(t_{n})$ is continuously differentiable. Furthermore, by the boundedness property \eqref{4.3} of an analytic semigroup $e^{-\tau A}$, it can be shown that 
\begin{align*}
	\left\| \int_0^{\tau} e^{-\tau A} (\tau-s) \tilde{w}(t_{n}+s) ds \right\| \le C \tau^2.
\end{align*}
In order to bound $d_{n+1}$, it is sufficient to bound $l_{n}^{\prime}(\sigma)$. The first derivative with respect to $\sigma$ is given by
\begin{align}\label{4.15}
	l_{n}^{\prime}(\sigma) = A e^{-(\tau -\sigma)A} \tilde{l}_{n}(\sigma) + e^{-(\tau-s)A} \tilde{l}_{n}^{\prime}(\sigma).
\end{align}
By the definition of $\theta$ given by \eqref{3.1}, it holds that $\| \theta\| \le K$, for all fixed integers $K >0$. Using the boundedness property of the analytic semigroup \eqref{4.3}, the second term of the right hand side of \eqref{4.15} is bounded. We observe that an additional $A$ on the first term on the right hand side of \eqref{4.15} can be compensated by applying the parabolic smoothing property \eqref{4.4} in the error recursion for the global error (see the proof of Theorem \ref{Th: 4.1} below). Thus, the adapted Lie splitting has a consistency error proportional to $\tau$. That is
\begin{align}\label{4.16}
	d_{n+1} = A \tilde{d}_{n+1}+ \mathcal{O}(\tau^2), \quad \tilde{d}_{n+1} = \mathcal{O}(\tau^2),
\end{align}
which completes the local error estimate for the adapted Lie splitting method applied to \eqref{1.1}.

\subsection{Global error}
The first-order convergence of the adapted Lie splitting applied to \eqref{1.1} is given in the next theorem. Before giving our main convergence analysis, we will make the following assumptions on the data of the problem \eqref{1.1}.
\begin{assumption}\label{as2}
	Let $D$ be a second-order strongly elliptic differential operator with smooth coefficients, $f$ and $b$ are twice continuously differentiable, and the convective coefficient $c(x)$ satisfies the condition \eqref{2.5}. Moreover, we also assume that the initial data $u_0$ are spatially smooth and satisfy the boundary conditions.
\end{assumption}
As a consequence, under the aforementioned assumptions, Theorem \ref{Th: 2.2} implies that the solution $u(t)$ of \eqref{1.1} is continuously differentiable. This allows us to state our main convergence theorem.
\begin{theorem} \label{Th: 4.1}(Convergence of the adapted Lie splitting).
	Let the Assumptions \ref{as1}-\ref{as2} be in charge. Then there exists a constant $\tau_0>0$ such that for all step sizes $0<\tau \le \tau_0$ and $t_{n}=n\tau$ we have that the adapted Lie splitting \eqref{3.8} applied to \eqref{1.1} satisfies the global error bound
	\begin{align}
		\Vert u_{n}-u(t_{n}) \Vert \le C\tau(1 + \left| \log \tau \right| ), \ \ 0\le n\tau \le T,
	\end{align}
	where the constant $C$ depends on $T$ but is independent of $\tau$ and $n$.
\end{theorem}
\begin{proof}
Let us denote by $e_{n}=u_{n}-u(t_{n})$ the global error. Thus, we have
\begin{align}\label{4.18}
	e_{n+1}=\mathcal{L}_{\tau}^{Ad}u_{n}-\mathcal{L}_{\tau}^{Ad} u(t_{n})+d_{n+1},
\end{align}	
where $d_{n+1}= \mathcal{L}_{\tau}^{Ad} u(t_{n})- u(t_{n+1})$ denotes the local error. 

By \eqref{4.8} and \eqref{4.11}, we have
\begin{align}\label{4.19}
	\mathcal{L}_{\tau}^{Ad}u_{n}-\mathcal{L}_{\tau}^{Ad} u(t_{n})
	= e^{-\tau A} e_{n} + 	\tau R(u_{n},u(t_{n}))+ \mathcal{O}(\tau^2),
\end{align}
where 
\begin{equation}\label{4.20}
	\begin{aligned}
		R(u_{n},u(t_{n}))=  e^{-\tau A} \Big( \nabla \cdot (\theta(\cdot) (u_{n}-u(t_{n})) ) \Big).
	\end{aligned}
\end{equation}
Inserting \eqref{4.19} and \eqref{4.20} into \eqref{4.18} yields
\begin{align}\label{4.21}
	e_{n+1}=d_{n+1}+\tau R(u_{n},u(t_{n})) + \mathcal{O}(\tau^2),
\end{align}
As $\| e_0\| =0$, solving the error recursion \eqref{4.21} we get
\begin{equation}
	\begin{aligned}
		e_{n}&= \sum_{k=0}^{n-1} e^{-(n-k-1)\tau A}(d_{k+1} + \mathcal{O}(\tau^2)) +  \tau \sum_{k=0}^{n-1} e^{-\left( n-k-1\right) \tau A}  R(u_{k},u(t_{k})).
	\end{aligned}
\end{equation}
The estimation of the local error $\Vert d_{k+1}\Vert$, given by \eqref{4.16} allows us to obtain
\begin{equation*}
	\begin{aligned}
		\Vert e_{n} \Vert&\le  \sum_{k=0}^{n-1}  \left\| e^{-(n-k-1)\tau A} (d_{k+1} + \mathcal{O}(\tau^2)) \right\| +  \tau \sum_{k=0}^{n-1} \left\| e^{-\left( n-k-1\right) \tau A}  R(u_{k},u(t_{k})) \right\| \\
		& \le \sum_{k=0}^{n-2} \left\| e^{-(n-k-1)\tau A} (A \mathcal{O}(\tau^2)+ \mathcal{O}(\tau^2)) \right\|  + \| d_{n} \| +C \tau \sum_{k=1}^{n-1} \left\|  e^{-(n-k)\tau A} \nabla \cdot (\theta(\cdot) e_{k})\right\|.
	\end{aligned}
\end{equation*}
Since $\| d_{n}\| \le C\tau$, using the parabolic smoothing property \eqref{4.4}, we arrive at
\begin{equation*}
	\begin{aligned}
		\| e_{n}\| & \le C\tau^2 \sum_{k=0}^{n-2} \left\| e^{-(n-k-1)\tau A} A\right\| + nC \tau^2 + C\tau \\
		& \qquad +C	\tau \sum_{k=1}^{n-1} \left\| e^{-(n-k)\tau A} A^{\frac{1}{2}} \right\| \| A^{-\frac{1}{2}} \nabla \| \| e_{k}\| \\
		&\le   C\tau^2 \sum_{k=1}^{n-1} \frac{1}{k\tau} + C\tau+ C\tau \sum_{k=1}^{n-1}  \frac{1}{\Big( (n-k) \tau \Big)^{\frac{1}{2}}}\Vert e_{k}\Vert   \\
		& \le C\tau(1+ \log \left| \tau\right| )  + C\tau \sum_{k=1}^{n-1}  t_{n-k}^{-\frac{1}{2}} \Vert e_{k}\Vert.
	\end{aligned}
\end{equation*}
This yields the desired bound by employing a discrete Gronwall lemma. This completes the proof.
\end{proof}

	
\section{Numerical experiments in two space dimension} 	
In this section, we present a series of numerical results for the convection-diffusion problem \eqref{1.1} on the domain $\Omega = [-0.5,1]^2$ , where the diffusion term is modeled by the Laplacian with the diffusion coefficient set to $0.01$ and the source term $f \equiv 0$. The unbounded convective term is set to $c(x, y) = 1/ \vert x \vert + 1/ \vert y \vert $ which belongs to the Lorentz spaces. The Laplacian is discretized by standard centered second-order finite differences, and the convection term is discretized by a first-order upwind scheme with $99$ grid points with an equidistant mesh. The reference solution is computed using the $\mathrm{ODE45}$ with absolute and relative tolerances are set to $10^{-12}$. In all our simulations, the scheme \eqref{1.4} is referred to as the classical Lie splitting, while the scheme proposed in Section~3.2 is referred to as the adapted Lie splitting.

\begin{figure}[ht]
	\centering
	\subfigure[Time-independent BCs \label{fig:subfig_a}]{\includegraphics[width=0.48\textwidth]{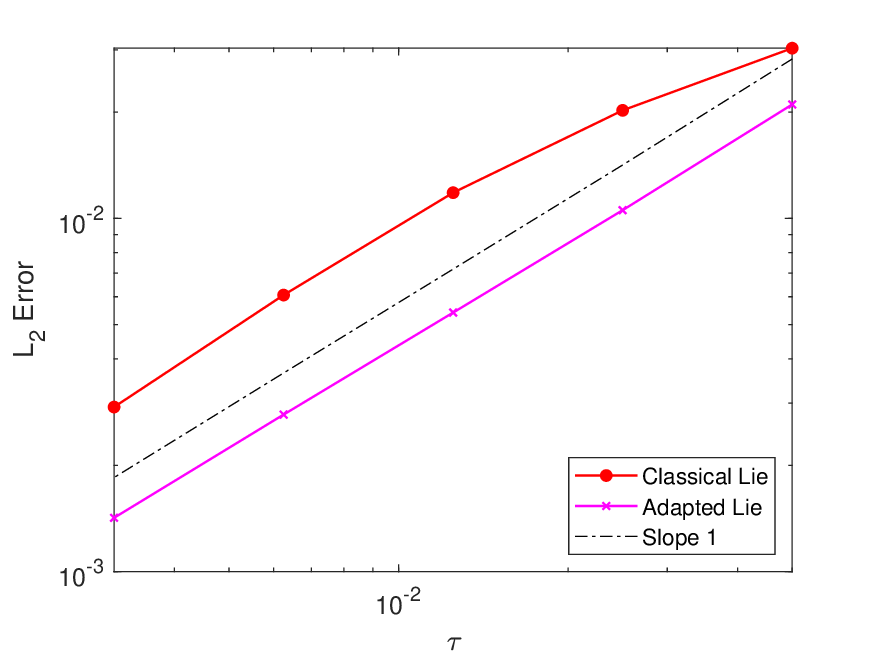}}
	\subfigure[Time-dependent BCs \label{fig:subfig_b}]{\includegraphics[width=0.48\textwidth]{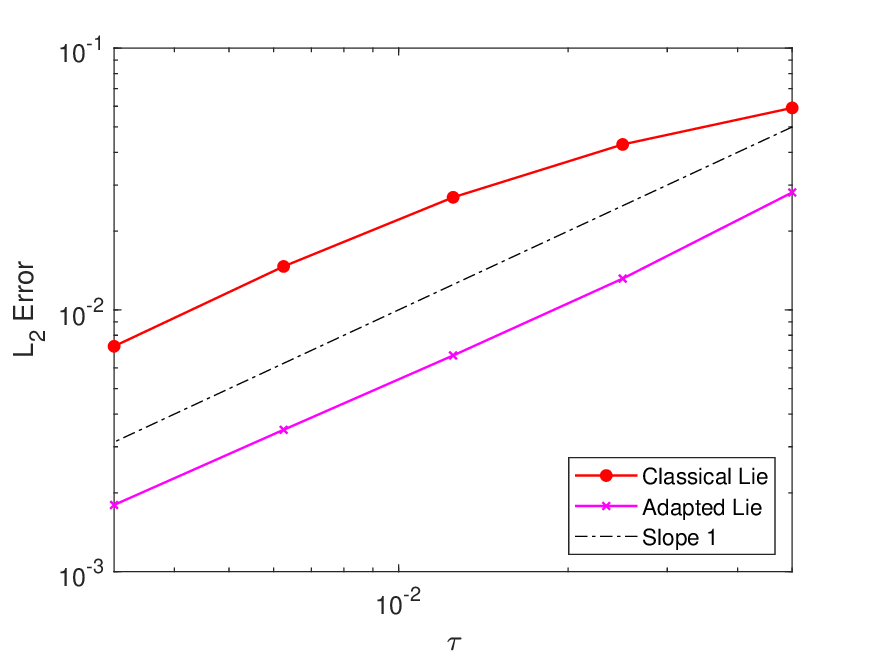}}	\\		

	\caption{ A comparison between the Lie splitting method and the newly proposed adapted Lie splitting method when applied to convection-diffusion equation with time-independent (Figure \ref{fig:subfig_a}) and time-dependent (Figure \ref{fig:subfig_b}) boundary conditions. The error in the discrete $L^2$ norm is computed at time $T=0.1$ and is plotted by comparing the numerical solution to a reference solution. A dash-dotted slope one has been added for comparison.}
	\label{fig1}
\end{figure}

\begin{exam} (time-independent boundary conditions) This problem considers the convection-diffusion equation with constant Dirichlet boundary conditions, which is given by $ u(t,x,-\frac{1}{2})=u(t,x,1) =u(t,-\frac{1}{2},y) = u(t,1,y)=1$. We have chosen the initial data $u(0,x, y)= 1+ \sin \left(   \frac{2}{3} \pi \left( x+\frac{1}{2}  \right) \right) \sin \left(   \frac{2}{3} \pi \left( y+\frac{1}{2}  \right) \right)$ that satisfies the presribed boundary conditions. Our numerical results are shown in Figure \ref{fig:subfig_a}.
\end{exam}	
\begin{exam}(time-dependent boundary conditions) In this example, we prescribe the convection-diffusion equation with time-dependent boundary conditions as follows: $u(t,-\frac{1}{2},y) = u(t,1,y) = 1 + \sin(5t)$ and $ u(t,x,-\frac{1}{2}) = u(t,x,1) = 1 + \sin(10t).$ As the initial data we have chosen $u(0,x, y)= 1+ \sin \left(   \frac{2}{3} \pi \left( x+\frac{1}{2}  \right) \right) \sin \left(   \frac{2}{3} \pi \left( y+\frac{1}{2}  \right) \right)$. The numerical results are displayed in the Figure \ref{fig:subfig_b}.
\end{exam}	
As expected, Figure \ref{fig1} illustrates the instability effect observed in the classical Lie splitting (see the red line), whereas the adapted Lie splitting (see the magnetic line) appears to be perfectly stable in both cases with time-independent and time-dependent boundary conditions. Furthermore, our numerical simulations have demonstrated that for the time-independent boundary conditions, the adapted Lie splitting method achieves an accuracy factor of approximately $2$ smaller than the classical Lie splitting method with the same time step size. In the case of time-dependent boundary conditions, the improved accuracy obtained by the adapted Lie splitting method is approximately a factor of $4$ compared to the classical Lie splitting method. The lower accuracy is maintained when a smaller time step size is employed, as the two methods have an order of convergence of one. It is evident that the adapted Lie splitting method offers advantages over the traditional Lie splitting method, with greatly improved accuracy. 

\begin{figure}[ht]
	\centering
	\subfigure[The mesh \label{fig2:subfig_a}]{\includegraphics[width=0.48\textwidth]{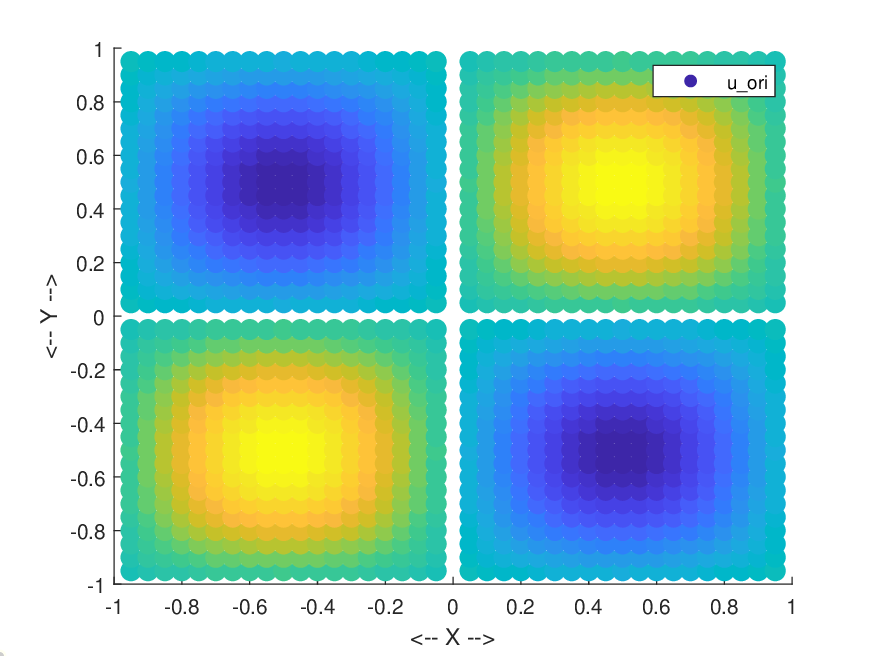}
	}
	\subfigure[The reference solution \label{fig2:subfig_b}]{\includegraphics[width=0.48\textwidth]{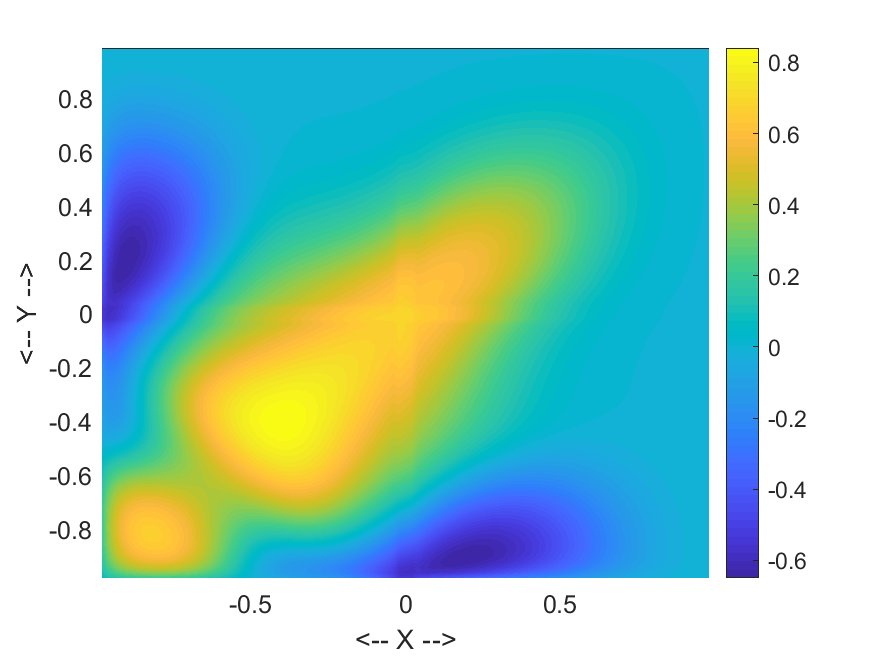}}	\\	
	
	\subfigure[The error  \label{fig2:subfig_c}]{\includegraphics[width=0.48\textwidth]{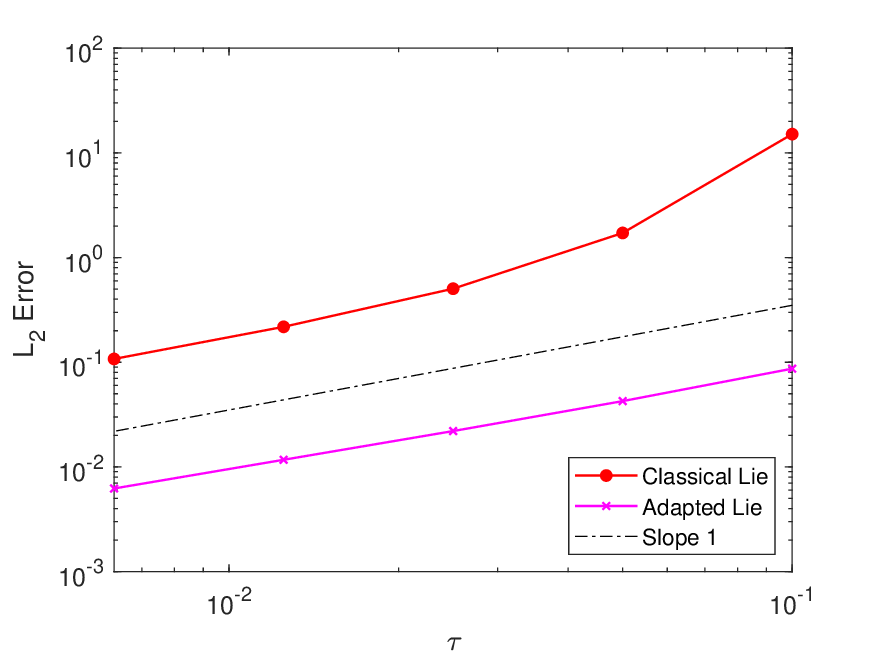}}				
	\caption{ A comparison between the Lie splitting method and the newly proposed adapted Lie splitting method when applied to the convection-diffusion problem \eqref{1.1} with homogeneous boundary conditions. The error at time $T=0.1$ is computed by comparing the numerical solution to the reference solution and is plotted as a function of the time step sizes, $\tau = \frac{1}{10},\frac{1}{20},\frac{1}{40},\frac{1}{80},\frac{1}{160}$. A dash-dotted slope one is drawn for comparison. The reference solution is shown on the right picture \ref{fig2:subfig_b}.}
	\label{fig2}
\end{figure}

	\begin{exam}(homogeneous boundary conditions) In this example, we consider the convection-diffusion problem on the domain $\Omega = [-1,1]^2$. The diffusion is modeled by the Laplacian with the confusion coefficient set to $0.1$ and the convection is also set to $c(x,y)=c(x,y) = 1/ \vert x \vert + 1/ \vert y \vert $. We impose the equation with homogeneous Dirichlet boundary conditions. The initial data $u_0(x,y) = \sin \left( \pi x \right) \sin \left( \pi y \right)$ is chosen to satisfy the boundary conditions. For the implementation, the equidistant mesh is used to discretize in space. Consequently, it is necessary to remove the element $0$ in both the $x$-axis and the $y$-axis to avoid singularity (see Figure \ref{fig2:subfig_a}). Our numerical results are shown in Figure \ref{fig2}.
\end{exam}	
It can be observed that the classical Lie splitting (the red line) is unstable as expected while the adapted Lis splitting (the magnetic line) appears to be perfectly stable (see Figure \ref{fig2}c). Furthermore, it can be seen that the improved accuracy obtained by the adapted Lie splitting method is approximately a factor of $18$ over the classical Lie splitting method. 
\section{Numerical experiments in three space dimension}
\begin{figure}{ht}
		\centering
			\includegraphics[width= 0.5\linewidth]{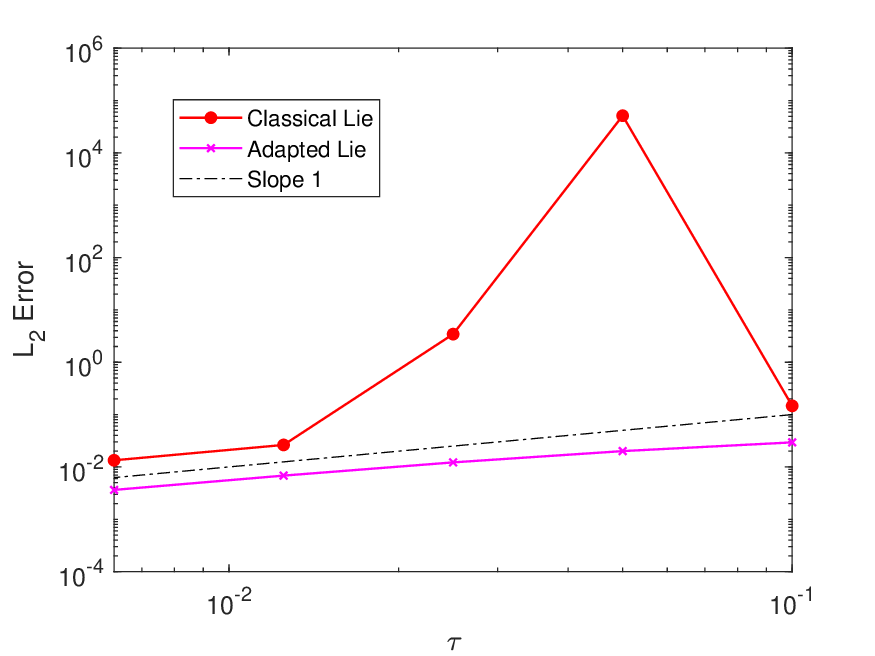}
		\caption{A comparison between the Lie splitting method and the newly proposed adapted Lie splitting method when applied to the convection-diffusion problem \eqref{1.1} with homogeneous boundary conditions in three dimension space. The absolute error in the discrete $L^2$ norm is computed at $T=0.1$ by comparing the numerical solution to a reference solution.}
		\label{fig3}
\end{figure}

In this section, we present numerical results for the convection-diffusion problem \eqref{1.1} for the unknown $u(t,x)$ in $[0,T]\times [-0.5,1]^3$ of the form:
\begin{equation*}\left\{
	\begin{aligned}
		& \partial_{t}u- \triangle u- \mathrm{ div} \left\lbrace \left[ \gamma \frac{x}{\vert x \vert^2} + \frac{1}{2-N+ \gamma} \Big( 1- \vert x \vert^{N-2-\gamma}\Big)x \right] u\right\rbrace  = - \mathrm{ div} \left(  \frac{x}{\vert x \vert^{N -\gamma}} \right) , \\
		& u|_{\partial \Omega}=0,\\
		& u (0)=0,
	\end{aligned}
	\right.
\end{equation*}
where $\frac{N}{2} -1 < \gamma \le N-1$ with $\gamma \ne N-2$. In the simulation, we have selected $N=3$ and $\gamma = 2$. The Laplacian is approximated by second-order finite differences methods, and the convection is discretized by a first-order upwind scheme with $19$ grid points with an equidistant mesh. The numerical results are presented in Figure \ref{fig3}, which confirms that the classical Lie splitting (the red line) exhibits significant instability effects, in contrast to the remarkably improved stable behavior of the adapted Lie splitting (the magnetic line), which allows for large time step sizes. Furthermore, the adapted Lie splitting method shows superior accuracy compared to the classical Lie splitting.

\section{Conclusion and outlook}
We have proposed a newly adapted splitting approach that has proven successful in overcoming the instability encountered in convection-diffusion problems with unbounded convective coefficients. Our numerical experiments have demonstrated that the adapted Lie splitting is more accurate compared to the standard Lie splitting method. In addition, our proposed approach is potentially applicable to large systems of convection-diffusion, with a wider range of applicability in real-world phenomena. Extending our proposed approach to second-order convergent Strang splitting, or even higher-order splitting methods, has the potential to be beneficial, but it introduces additional complexities due to the unbounded property of the convective term. Such extensions may necessitate higher regularity assumptions on the solution. Furthermore, incompatibility effects may arise within the internal step of the Strang splitting scheme, resulting in order reduction (see, e.g., \cite{doi:10.1137/19M1257081,doi:10.1137/16M1056250}). Nevertheless, despite the aforementioned complexities, the potential for enhanced convergence properties and computational efficiency of adapted splitting methods remains a compelling topic for future research.

\nocite{*}
\bibliographystyle{siam}
\bibliography{references_DHO}

\end{document}